\documentclass{birkmult}

\usepackage{amssymb}
\usepackage{amsmath}
\usepackage{latexsym}

\usepackage{amscd}
\usepackage{amsfonts}
\usepackage[all]{xy}
\usepackage{mdwtab}
\usepackage{float}
 \newtheorem{Thm}{Theorem}[section]
 \newtheorem{Cor}[Thm]{Corollary}
 \newtheorem{Lem}[Thm]{Lemma}
 \newtheorem{Pro}[Thm]{Proposition}
 
 \newtheorem{Fac}[Thm]{Fact}

\theoremstyle{remark}
\newtheorem*{War}{Warning}
\theoremstyle{definition}
                             \newtheorem{Rem}[Thm]{Remark}
                             
                             \newtheorem{Def}[Thm]{Definition}

                             \newtheorem{Ex}[Thm]{Example}
                             \newtheorem{Que}[Thm]{Question}

\newtheorem*{ConNot}{Conventions and Notations}

\renewcommand{\thetable}{\theThm}
\newcommand{\longrightsquigarrow}
           {\smallsmile\!\smallsmile\!\smallsmile\!\smallsmile\!\smallsmile\!\smallsmile\!\smallsmile\!\smallsmile\!\rightarrow}
\newcommand{\ZZ}{\mathbb {Z}}

\newcommand{\RR}{\mathbb R}
\newcommand{\CC}{\mathbb C}
\newcommand{\RP}{\mathbb{RP}}
\def\C{\mathrm {C}}
\def\N{\mathrm {N}}
\renewcommand{\rho}{\varrho}
\newcommand{\eps}{\varepsilon}
\renewcommand{\d}{\partial}
\renewcommand{\:}{\colon\,}

\newcommand{\id}{\mathrm {id}}
\newcommand{\fr}{\mathrm {fr}}
\newcommand{\inv}{\mathrm {inv}}
\newcommand{\forg}{\mathrm {forg}}
\newcommand{\coll}{\mathrm {coll}}
\newcommand{\im} {{\mathrm {im}}}
\newcommand{\Int}[1]{{\stackrel{{\scriptscriptstyle\circ}}{#1}}}

\newcommand{\Proof} {\noindent{\itshape Proof. \ }}
\renewcommand{\proof}[1] {\noindent{\itshape Proof {#1}. \ }}
\newcommand{\pproof}[1] {\noindent{{\itshape Proof }({#1}). \ }}

\newcommand{\Bigdownarrow}{\Downarrow}

\begin{document}
%
%
%
%
%
%
%
%
%
\title[Selfcoincidences and roots in Nielsen theory]
{Selfcoincidences and roots in Nielsen theory}
\author[Ulrich Koschorke]{Ulrich Koschorke}

\address{%
Fachbereich Mathematik\\
Universit\"at Siegen\\
57068 Siegen - Germany}

\email{koschorke@mathematik.uni-siegen.de}

\thanks{Supported in part by DAAD (Germany) and CAPES (Brazil).}
\subjclass{Primary 54H25, 55M20; Secondary 55N22, 55P35}

\keywords{Coincidence invariants; Nielsen number; looseness
obstruction; Hopf invariant.}



\begin{abstract}
Given two maps $f_1$ and $f_2$ from the sphere $S^m$ to an
$n$-manifold $N$, when are they loose, i.e.\ when can they be
deformed away from one another? We study the geometry of their
(generic) coincidence locus and its Nielsen decomposition. On the
one hand the resulting bordism class of coincidence data and the
corresponding Nielsen numbers are strong looseness obstructions. On
the other hand the values which these invariants may possibly assume
turn out to satisfy severe restrictions, e.g.\ the Nielsen numbers
can only take the values $0,1$ or the cardinality of the fundamental
group of $N$. In order to show this we compare different Nielsen
classes in the root case (where $f_1$ or $f_2$ is constant) and we
use the fact that all but possibly one Nielsen class are inessential
in the selfcoincidence case (where $f_1=f_2$). Also we deduce strong
vanishing results.
 \end{abstract}

\maketitle

\section{Introduction} \label{sec:intro}

This work is dedicated to Albrecht Dold whose enormous influence on
the development of topology continues to this very day. In fact, it
was his stimulating paper~\cite {DG} (written jointly with Daciberg
Gon\c calves) which got me interested in coincidences of maps and
which inspired me to introduce a related obstruction theory
(cf.~\cite {K2}--\cite {K7}). Here we apply it to the special case
where the domain is a sphere.

Throughout this paper let $f_1,f_2\:S^m\to N$ be two (continuous)
maps into an $n$-dimensional smooth connected manifold $N$ without
boundary. We will also assume that $m,n\ge2$ (the remaining cases
being well understood).

We want to make the coincidence locus
$$ \refstepcounter {Thm} \label {eq:1.1}
\C(f_1,f_2) = \{x\in S^m\mid f_1(x)=f_2(x)\}
 \leqno (\theThm)
$$
as small as possible (in some sense) by deforming $f_1$ and $f_2$.

\begin {Def}  \label {def:1.2}
(compare~\cite {DG}, p.~293--296 and~\cite {GR2}, \S1, where a
different terminology is used for the same concepts; see also~\cite
{K7}, 5.3).

(i) The pair $(f_1,f_2)$ is called {\it loose} if there exist maps
$f_i'$ homotopic to $f_i$, $i=1,2$, whose coincidence locus
$\C(f_1',f_2')$ is empty (i.e.\ ``$f_1$ and $f_2$ can be deformed
away from one another'').

(ii) In the special case $f_1=f_2=:f$ we have a refined notion: the
pair $(f,f)$ is {\it loose by small deformation} if for every metric
on $N$ and for every $\eps>0$ there exists an $\eps$-approximation
$\bar f$ of $f$ such that $\bar f(x)\ne f(x)$ for all $x\in S^m$.
 \end{Def}

 \begin{Pro} \label {pro:1.3}
All pairs of maps $f_1,f_2\:S^m\to N$ are loose if at least one of
the following conditions hold:

{\rm (i)} $m<n$; or

{\rm (ii)} $N$ is not compact; or

{\rm (iii)} the fundamental group $\pi_1(N)$ is not finite; or

{\rm (iv)} $N$ is the total space of a Serre fibration with a section
and with strictly positive dimensions of the fiber and base (e.g.\ if
$N$ is the product of such spaces).
\end{Pro}

 \begin{Pro} \label {pro:1.4}
Assume that

{\rm (v)} $N$ is noncompact or has trivial Euler characteristic
$\chi(N)$ (this holds e.g.\ when $n$ is odd), or

{\rm (vi)} $\pi_{m-1}(S^{n-1})=0$.

Then for all maps $f\:S^m\to N$ the pair $(f,f)$ is loose by small
deformation.
\end{Pro}

\begin{Ex} [surfaces] \label {ex:1.5}
Let $n=2$ and $m\ge1$. Then all such pairs $(f,f)$ are loose by small
deformation except when $m=2$ and $N=S^2$ or $\RP(2)$.
\end{Ex}

It is not hard to verify the looseness claims in~\ref {pro:1.3}, \ref
{pro:1.4} and~\ref {ex:1.5} by elementary considerations (see
section~\ref {sec:2} below). However, in general a deeper analysis is
needed and the following approach has proved to be fruitful.

For every pair of (base point preserving) maps $f_1,f_2\:S^m\to N$ we
defined (in~\cite {K6}) an invariant
$$ \refstepcounter {Thm} \label {eq:1.6}
\omega^\#(f_1,f_2) \in \pi_m\left(S^n\wedge(\Omega N)^+\right)
 \leqno (\theThm)
$$
which reflects the geometry of a (generic) coincidence submanifold
$C\subset S^m$, its normal bundle as well as certain path space data
(for details see~\cite {K6} or section~\ref {sec:3} below). The
invariant $\omega^\#(f_1,f_2)$ depends only on the (base point
preserving) homotopy classes of $f_1,f_2$ and must vanish if
$(f_1,f_2)$ is loose (and, in particular, in all cases listed in
proposition~\ref {pro:1.3}). The converse holds in a ``stable''
dimension range.

\begin{Thm} \label {thm:1.7}
{\rm (compare~\cite {K3}, 1.10 and~\cite{K2}, 2.2). } Assume
$m<2n-2$. Then a pair $(f_1,f_2)$ is loose precisely if
$\omega^\#(f_1,f_2)=0$. In the special case $f_1=f_2=:f$ the pair
$(f,f)$ is loose by small deformation precisely if
$\omega^\#(f,f)=0$.
 \end{Thm}

Our looseness obstruction is always compatible with additions in
homotopy groups:
$$ \refstepcounter {Thm} \label {eq:1.8}
\omega^\#(f_1+f_1',f_2+f_2') = \omega^\#(f_1,f_2) +
\omega^\#(f_1',f_2')
 \leqno (\theThm)
$$
for all $[f_i],[f_i']\in \pi_m(N)$, $i=1,2$.

There is also a canonical involution \ $\inv$ \ on the group
$\pi_m(S^n\wedge(\Omega N)^+)$. It plays a r\^ole e.g.\ in the
symmetry relation
$$ \refstepcounter {Thm} \label {eq:1.9}
\omega^\#(f_2,f_1) = \inv\left(\omega^\#(f_1,f_2)\right).
 \leqno (\theThm)
$$

\medskip
Two special settings are very central in coincidence theory.

\bigskip
{\bf I. The root case: $f_2=y_0$ } (where the fixed value $y_0\in N$
is given).

Here our invariant yields the degree homomorphism
$$ \refstepcounter {Thm} \label {eq:1.10}
\deg^\#:=\omega^\#(-,y_0)\ \:\ \pi_m(N)\to\pi_m(S^n\wedge(\Omega
N)^+).
 \leqno (\theThm)
$$
(For a purely homotopy theoretical interpretation in terms of an
enriched Hopf-Ganea invariant see theorem~7.2 in~\cite {K6}.)

Clearly $\deg^\#$ vanishes on $i_*(\pi_m(N\setminus\{*\}))$ where
$i\:N\setminus\{*\}\hookrightarrow N$ denotes the inclusion of the
complement of a point $*$ in $N$. It turns out that the resulting
sequence
$$ \refstepcounter {Thm} \label {eq:1.11}
 \pi_m(N\setminus\{*\}) \stackrel {i_*}{\longrightarrow}
 \pi_m(N) \stackrel {\deg^\#}{\longrightarrow}
 \pi_m(S^n\wedge(\Omega N)^+)
 \leqno (\theThm)
$$
is very often exact, e.g.\ when $m\le2n-3$ or $n\le2$ or $N$ is a
sphere or a (real, complex or quaternionic) projective space of
arbitrary dimension (cf.~\cite {K6}, (6.5)). Thus in these cases
$\deg^\#(f)$ is the {\it complete} looseness obstruction for the pair
$(f,y_0)$.

\begin {Que} \label {que:1.12}
Is the sequence~(\ref {eq:1.11}) {\it always} exact?
\end{Que}

The degree $\deg^\#$ can be considered to be our basic coincidence
invariant since (in view of~(\ref {eq:1.8}) and~(\ref {eq:1.9}))
$$ \refstepcounter {Thm} \label {eq:1.13}
 \begin {array} {lcl}
 \omega^\#(f_1,f_2)&=&\omega^\#(f_1,y_0) + \omega^\#(y_0,f_2) \\
                   &=&\deg^\#(f_1) + \inv(\deg^\#(f_2))
 \end{array}
 \leqno (\theThm)
$$
for all $[f_1],[f_2]\in\pi_m(N,y_0)$.

\bigskip
{\bf II. The selfcoincidence case: $f_1=f_2=:f$. }

Consider the bundle $ST(N)$ of unit tangent vectors over $N$ (with
respect to some metric) and the resulting exact (horizontal) homotopy
sequence as well as the Freudenthal suspension homomorphism $E$:
$$ \refstepcounter {Thm} \label {eq:1.14}
\xymatrix{
 \ldots \ar[r] & \pi_m(ST(N)) \ar[r] & \pi_m(N) \ar[r]^{\!\!\!\!\!\!\d} & \pi_{m-1}(S^{n-1}) \ar[r] \ar[d]^E & \ldots
 \\
               &                     &                    & \pi_m(S^n)                         & .
}
 \leqno (\theThm)
$$

\begin {Thm} \label {thm:1.15}
Given $[f]\in\pi_m(N)$, we have the following logical implications:

\medskip
{\rm (i)} $\d([f])\in\pi_{m-1}(S^{n-1})$ vanishes;

$\phantom{j}\!$$\Updownarrow$

{\rm (ii)} $(f,f)$ is loose by small deformation;

\medskip
$\phantom{i}$$\Bigdownarrow$ \quad $\left(\Updownarrow \mbox { if } N=\RP(n)\right)$

\medskip
{\rm (iii)} $(f,f)$ is loose (by any deformation);

\medskip
$\phantom{i}$$\Bigdownarrow$ \quad $\left(\Updownarrow \mbox { if } N=S^n\right)$

\medskip
{\rm (iv)} $\omega^\#(f,f)=0$;

$\phantom{i}$$\Updownarrow$

{\rm (v)} $E(\d([f])=0$.
\end{Thm}

Thus $\omega^\#(f,f)$ is just ``one desuspension short'' of being the
{\it complete} looseness obstruction.

The equivalence of~(i) and~(ii) was already noted by Dold and Gon\c
calves in~\cite {DG}.

Observe also that all conditions (i)--(v) in~\ref {thm:1.15} {\it
except} (iii) are compatible with covering projections $p\:\tilde
N\to N$ (compare~\cite {K7},1.22).

\begin{Cor} \label {cor:1.16}
The conditions {\rm (i)}$,\ldots,${\rm (v)} in~$\ref {thm:1.15}$ are
all equivalent if the suspension homomorphism $E$, when restricted to
$\d(\pi_m(N))$ (cf.~$(\ref {eq:1.14})$), is injective and, in
particular, if $m\le n+3$ or if $m=n+4\ne10$ or in the stable
dimensional range $m\le 2n-3$.
\end{Cor}

Indeed, in these three dimension settings $E$ is injective whenever
$n\equiv 0(2)$. Clearly conditions (i)--(v) in~\ref {thm:1.15} are
automatically satisfied under the assumptions of proposition~\ref
{pro:1.4}.

\begin{Ex} [$N=\RP(n)$ or $S^n$] \label{ex:1.17}
If $m\le n+4$, then the five conditions in~\ref {thm:1.15} are
equivalent for all maps $f\:S^m\to\RP(n)$ and $\tilde f\:S^m\to S^n$
(even in the exceptional case $m=n+4=10$ since $\pi_{10}(S^6)=0$).

However, this is no longer true for $m=n+5=11$. Indeed, according
to~\cite {T} and~\cite{P} we have in~(\ref {eq:1.14})
 $$
 \frac12H \: \pi_{11} (S^6) \stackrel {\cong}{\longrightarrow} \ZZ;
 \quad \pi_{10} (S^5) \cong \ZZ_2;
 \quad \pi_{10} (V_{7,2}) =0
 $$
where $H$ denotes the Hopf invariant. Thus $\d$ is onto, but $E$ and
hence $E\circ\d$ is trivial here. Therefore, given any map
$f\:S^{11}\to\RP(6)$ and a lifting $\tilde f\:S^{11}\to S^6$ of it,
the invariants $\omega^\#(f,f)$ and $\omega^\#(\tilde f,\tilde f)$
vanish and the pair $(\tilde f,\tilde f)$ is loose (see~\cite{K6},
1.12; compare also~\cite {GW}). But $(f,f)$ is loose and,
equivalently, $(\tilde f,\tilde f)$ is loose by small deformation
only (and precisely) if $[f]\in2\cdot\pi_{11}(\RP(6))$ or,
equivalently, if $H(\tilde f)\equiv0(4)$ (compare~\cite {GR1}
and~\cite{K7}, 1.19; the delicate difference between conditions (ii)
and (iii) in~\ref {thm:1.15} is further illustrated e.g.\
in~\cite{GR2} and~\cite{K7}, 1.22).

For simple examples of nontrivial $\omega^\#$-values consider the
case $m=n$. If a map $\tilde f\:S^n\to S^n$ has (standard) mapping
degree $d$, then
 $$
\pm\omega^\#(\tilde f,\tilde f) = E\circ\d(\tilde f) =
\left(1+(-1)^n\right)d
 $$
in $\pi_n(S^n\wedge(\Omega S^n)^+)\cong\pi_n(S^n)\cong\ZZ$. Assume
that $n$ is even. Then $(\tilde f,\tilde f)$ is loose if and only if
$\tilde f$ is nullhomotopic; the same holds for maps from $S^n$ into
$\RP(n)$. This shows that the exception made in example~\ref {ex:1.5}
is indeed necessary.
 \qed
\end{Ex}

\medskip
In previous work (cf.~\cite{K6}) we simplified our
$\omega^\#$-invariant by {\it assuming} that $f_1$ or $f_2$, say
$f_2$, is ``not coincidence producing'' (cf.~\cite{BS}), i.e.\
$$ \refstepcounter {Thm} \label {eq:1.18}
 (\bar f_2,f_2) \mbox{ is loose for {\it some} map }\bar f_2\:S^m\to N.
 \leqno (\theThm)
$$
Then $(f_1,f_2)$ and $(f_1-\bar f_2,(f_2-f_2)\sim y_0)$ have a
similar coincidence behaviour and
$$ \refstepcounter {Thm} \label {eq:1.19}
 \omega^\#(f_1,f_2) = \deg^\#(f_1-\bar f_2).
 \leqno (\theThm)
$$

In this paper we start from the decomposition (up to homotopy)
 $$
 (f_1,f_2) \sim (f_1-f_2,y_0) + (f_2,f_2)
 $$
which is always available without any assumption. In view of~(\ref
{eq:1.8}) this implies the basic equation
$$ \refstepcounter {Thm} \label {eq:1.20}
 \omega^\#(f_1,f_2) = \deg^\#(f_1-f_2)+\omega^\#(f_2,f_2).
 \leqno (\theThm)
$$
Vanishing results concerning the second (``selfcoincidence'') summand
are not only interesting in view of theorem~\ref {thm:1.15}, but they
also allow us to reduce our analysis to studying the degree
homomorphism $\deg^\#$ or, equivalently, Hopf-Ganea invariants
(cf.~\cite{K6}, 7.2).

It turns out that the cardinality $\#\pi_1(N)\in\ZZ\cup\{\infty\}$ of
the fundamental group of $N$ plays a crucial r\^ole.

\begin{Thm} \label {thm:1.21}
Assume that

{\rm (i)} $\#\pi_1(N)>2$ and $N$ is orientable or not; or

{\rm (ii)} $\#\pi_1(N)\ge2$ and $N$ is orientable.

Then $\omega^\#(f,f)=0$ for all $f\:S^m\to N$.

In particular
 $$
 \omega^\#(f_1,f_2)=\deg^\#(f_1-f_2)
 $$
for all $[f_1],[f_2]\in\pi_m(N)$.

Also the set of all possible values of our $\omega^\#$-invariant is
restricted by the relation
 $$
 \{\omega^\#(f_1,f_2)\}\ =\ \im(\deg^\#)\ \subset\ \ker(\id + \inv)
 $$
where \ $\id=$ identity and \ $\inv$ \ denotes the canonical
involution of $\pi_m(S^n\wedge(\Omega N)^+)$ (cf.~$(\ref {eq:1.9})$
and~$(\ref {eq:1.13})$).

If in addition $m<2n-2$ or $m\le n+3$, then the pair $(f,f)$ is loose
by small deformation for every map $f\:S^m\to N$.
 \end{Thm}

On the one hand, the proof (given in section~\ref {sec:4} below)
compares the contributions of the pathcomponents of the loop space
$\Omega N$ to $\omega^\#(f,f)$; on the other hand, it uses the fact
that at most one such pathcomponent can contribute nontrivially to
our selfcoincidence invariant.

\begin{Ex} \label{ex:1.22}
As an application consider the case where $m=n$ and $N$ is an
orientable $n$-manifold with $\pi_1(N)=\ZZ_2$. Then
 $$
\inv=(-1)^n\id \quad \mbox{on} \quad \pi_n(S^n\wedge(\Omega
N)^+)\cong\ZZ\oplus\ZZ.
 $$
If $n$ is even then \ $\ker(\id+\inv)=0$; therefore each pair of maps
$f_1,f_2\:S^m\to N$ is loose (since $\omega^\#(f_1,f_2)=0$).
\end{Ex}

\begin{Ex} [real Grassmann manifolds] \label{ex:1.23}
Consider the manifold $G_{r,k}$ (or $\tilde G_{r,k}$, resp.) of
nonoriented (or oriented, resp.) $k$-planes in $\RR^r$, $0<k<r$. If
$r$ is even, then $G_{r,k}$ is orientable, $\pi_1(G_{r,k})\ne\{0\}$
and hence $\omega^\#(f,f)=0$ for all maps from $S^m$, $m\ge1$, to
$G_{r,k}$ or $\tilde G_{r,k}$. This does not seem to follow
from~\ref {pro:1.3} or~\ref {pro:1.4}; indeed, the Euler
characteristic of $G_{r,k}$ is strictly positive if both $r$ and $k$
are even (see~\ref {fac:4.3} below; compare also our example~\ref
{ex:2.2}).
\end{Ex}

A further vanishing result (quoted from~\cite{K6},1.9; compare
also~\cite{K7}, 5.4) is of interest here.

\begin{Pro} \label {pro:1.24}
Assume $\pi_1(N)\ne0$. Given $[f]\in\pi_m(N)$, $m\ge2$, if
$j_*\circ\d([f])$ vanishes then so does $\omega^\#(f,f)$.
\end{Pro}

Here $j\:S^{n-1}\hookrightarrow N\setminus\{*\}$ denotes the
inclusion of the boundary sphere of a small $n$-ball in $N$ centered
at a point $*\in N$. According to~\cite{K7}, 5.4, the condition
$j_*\circ\d([f])=0$ just means that $f$ is not coincidence producing
(cf.~(\ref {eq:1.18})).
 \qed

\medskip
Let us summarize: {\it although our selfcoincidence invariant
$\omega^\#(f,f)$ is ``only one desuspension short'' of being a
complete looseness obstruction, it vanishes for a great number of
maps $f$ defined on spheres} (in contrast, in~\cite{K2}
$\omega^\#(f,f)$ was shown to be highly nontrivial in many examples
where the domain of $f$ is a general closed manifold $M$).

\begin{Cor} \label {cor:1.25}
Assume $\pi_1(N)\ne0$ and $n\ge1$. If $\omega^\#(f,f)\ne0$ for a map
$f\:S^m\to N$ then the following restrictions must all be satisfied:

{\rm a)} $n$ is even and $m\ge n\ge 4$, or else $m=2$ and $N=\RP(2)$;
and

{\rm b)} $N$ is closed and nonorientable, \ $\pi_1(N)=\ZZ_2$, \
$\chi(N)\ne0$; \ moreover the homomorphism
$i_*\:\pi_m(N\setminus\{*\},y_0)\to\pi_m(N,y_0)$ (induced by the
inclusion of $N$, punctured at some point $*\ne y_0$) is not onto;
furthermore, $N$ does not allow a free smooth action by a nontrivial
finite group; and

{\rm c)} $E\circ\d([f])\ne0$ and $j_*\circ\d([f])\ne0$
(compare~$(\ref {eq:1.14})$ and~$\ref {thm:1.21}$).
\end{Cor}

Thus an obvious place to look for nontrivial values of
$\omega^\#(f,f)$ are even-dimensional real projective spaces (cf.\
already example~\ref {ex:1.17}).

\begin{Ex} \label {ex:1.26}
If $n=4,8,12,14,16$ or $20$, then there exist infinitely many
homotopy classes $[f]\in\pi_{2n-1}(\RP(n))$ such that
$\omega^\#(f,f)\ne0$.

We see this with the help of the weaker invariant
$\omega(f,f)\in\pi_{m-n}^S$ (cf.~(\ref {eq:3.15}) below
and~\cite{K2}, 2.1), which stabilizes and simplifies $\omega^\#(f,f)$
and is more easily computable.

However, for many maps $f,f_1,f_2$ from $S^m$ into a nonorientable
$n$-manifold $N$ this approach does not yield nontrivial
$\omega$-values. Indeed, $\omega(f,f)=0$ whenever
$2\cdot\pi_{m-n}^S=0$, e.g.\ when $m-n=1,2,4,5,6,8,9,12,14,16$ or
$17$ (cf.~\cite{T}). Moreover, for all $m$ and $n$
$$ \refstepcounter {Thm} \label {eq:1.27}
 (\chi(N)-1)\ \omega(f_1,f_2) = 0
 \leqno (\theThm)
$$
where $\chi(N)$ stands for the Euler characteristic of $N$.
 \end{Ex}

\begin{Ex} [\rm $\mathbf{N=G_{5,2}}$; compare~\ref {ex:1.23}] \label {ex:1.28}
The invariant $\omega(f_1,f_2)$ vanishes for all maps
$f_1,f_2\:S^m\to G_{5,2}$, $m\ge1$. In particular, the induced
homomorphism
 $$
 \coll_*\:\pi_m(G_{5,2})\to\pi_m(S^6)
 $$
is trivial for $m\le10$ where \ $\coll$ \ denotes the degree one map
which collapses all but an open topdimensional ball into a point.
(Note that e.g.\ for $m=6,7$ or $8$, resp.,
$\pi_m(G_{5,2})\cong\pi_m(V_{5,2})$ is isomorphic to $\ZZ_2$,
$\ZZ\oplus\ZZ_2$ and $\ZZ_2$, resp.; cf.~\cite{P}).
\qed
\end{Ex}

In section~\ref {sec:4} below we will deduce the vanishing
theorem~\ref {thm:1.21} and relations such as~(\ref {eq:1.27})
and~\ref{ex:1.28} from
a careful analysis of the root case in section~\ref {sec:3}.

While our looseness obstructions lie in complicated groups which are
usually hard to compute, they give rise to simple numerical
invariants (defined in section~\ref {sec:5} below). These generalize
the Nielsen numbers which have played such a central r\^ole in
topological fixed point theory (cf.\ e.g.~\cite{B} and~\cite{K6},
1.12~(iv)). In analogy to the classical procedure, our Nielsen
numbers $\N^\#(f_1,f_2)\ge\N(f_1,f_2)$ count those ``Reidemeister
classes'' $A\in\pi_0(\Omega N)=\pi_1(N)$ which make an essential
contribution to $\omega^\#(f_1,f_2)$ and $\tilde\omega(f_1,f_2)$,
resp.

\begin{Thm} \label {thm:1.29}
{\rm (cf.~\cite{K6}, 1.2). } For every pair $(f_1',f_2')$ of maps
homotopic to $(f_1,f_2)$ the number of pathcomponents of the
coincidence subspace $\C(f_1',f_2')\subset S^m$ is at least
$\N^\#(f_1,f_2)$.
\end{Thm}
There are quite a few examples where $\N^\#(f_1,f_2)$ is the best
possible such lower bound (see e.g.\ the ``Wecken theorems''
in~\cite{K6}, 1.3, 1.12, 1.13, 1.14$,\ldots$).

\begin{Thm} \label {thm:1.30}
For all maps $f_1,f_2\:S^m\to N$ we have:
 $$
\begin{array}{ccc}
 \N^\#(f_1,f_2)=0 & \mbox{if and only if} & \omega^\#(f_1,f_2)=0\ ; \\
 \N(f_1,f_2)=0    & \mbox{if and only if} & \tilde\omega(f_1,f_2)=0.
 \end{array}
 $$
\end{Thm}

We call this the norm property of our Nielsen numbers: they decide
whether elements in the image group of our $\omega$-invariants are
zero, just like norms of vectors do. (Recall, in addition, analogues
of the triangular inequality, cf.~\cite{K6}, 6.2.)

By definition $\N^\#(f_1,f_2)$ is among the integers between $0$ and
$\#\pi_1(N)$. But few of these actually occur as Nielsen numbers.

\begin{Thm} \label{thm:1.31}
Let $k\in\ZZ\cup\{\infty\}$ be the number of elements in $\pi_1(N)$.
Then for each pair of maps $f_1,f_2\:S^m\to N$ the Nielsen numbers
$\N^\#(f_1,f_2)$ and $\N(f_1,f_2)$ may assume only the two values $0$
or $k$ or, if $k=2$ and $N$ is an even-dimensional, closed,
nonorientable manifold with nontrivial Euler characteristic, also $1$
as a third possible value.
\end{Thm}

Here again (as in~\cite{K7}) the case $N=\RP(n)$ turns out to be
particularly interesting. E.g.\ assume that $m=n$ is even and let
$p\:S^n\to\RP(n)$ be the canonical covering map. Then
$\N^\#(f_1,f_2)=\N(f_1,f_2)$ equals $0,1$ and $2$, resp., when
$(f_1,f_2)=(y_0,y_0)$, $(p,p)$ and $(p,y_0)$, resp.\ (cf.\ the end
of section~\ref {sec:4} below).

\begin{Rem} \label {rem:1.32}
It will be convenient to consider base point preserving maps and
homotopies in sections~\ref {sec:3}--\ref {sec:5} below. Recall,
however, that looseness and Nielsen numbers depend only on free
homotopy classes and so does the vanishing of our $\omega$-invariants
(see e.g.~\cite{K6}, 1.2, 2.1 and the appendix there).
\end{Rem}

\begin{Rem} \label {rem:1.33}
Our approach can also be extended to the setting of fibre preserving
maps between smooth fibrations. The resulting coincidence invariants
are closely related to A.~Dold's fixed point index which was defined
and studied in~\cite{D}.
\end{Rem}

\begin{ConNot}
Throughout this paper $N$ denotes a smooth connected manifold without
boundary (Hausdorff and having a countable basis). Our notation will
often not distinguish between a constant map and its value.
\end{ConNot}

\subsection*{Acknowledgment}
It is a pleasure to thank D.~Gon\c calves and E.~Kudryavtseva for
stimulating discussions.

\section{Looseness} \label{sec:2}

In this section we use rather elementary techniques to establish the
looseness results in propositions~\ref {pro:1.3} and~\ref {pro:1.4}
as well as example~\ref {ex:1.5}.

\begin{Lem} \label {lem:2.1}
Let $y_1\ne y_2$ be different points in $N$ and assume that the
homomorphism \ $i_*\:\pi_m(N\setminus\{y_1\},y_2)\to\pi_m(N,y_2)$
induced by the inclusion map is onto. Then for all maps
$f_1,f_2\:S^m\to N$ the pair $(f_1,f_2)$ is loose.
\end{Lem}

\Proof
 Let $S_+^m$ and $S_-^m$ denote the hemispheres defined by $x_1\ge0$
and $x_1\le0$, resp., $x=(x_1,\dots,x_{m+1})\in S^m$. Then $f_2$ is
homotopic to a map $f_2'$ such that $f_2'(S_+^m)\subset
N\setminus\{y_1\}$ and $f_2'|_{S_-^m}\equiv y_2$. Similarly, $f_1\sim
f_1'$ such that $f_1'|_{S_+^m}\equiv y_1$ and $f_1'(S_-^m)\subset
N\setminus\{y_2\}$. Clearly the pair $(f_1',f_2')$ is coincidence
free.
 \qed

\medskip
Proposition~\ref {pro:1.3} follows as a corollary since each of the
conditions~(i) through~(iv) imply that \ $i_*$ is onto. If $m<n$,
this is seen by making a map $f\:S^m\to N$ transverse to $\{y_1\}$.
If $N$ is not compact, deform $f$ by an isotopy along a smoothly
embedded path which starts in $N\setminus f(S^m)$, ends in $y_1$ and
avoids $y_2$.

If the universal covering space $p\:\tilde N\to N$ has infinitely
many layers, consider a lifting $\tilde f\:S^m\to\tilde N$ of $f$.
Apply an isotopy along suitable disjoint paths in $\tilde N$ which
start in $\tilde N\setminus\tilde f(S^m)$ and end in the finitely
many points of $\tilde f(S^m)\cap p^{-1}(\{y_1\})$. The corresponding
homotopy in $N$ deforms $f$ into $N\setminus\{y_1\}$.

Finally, in case~(iv) of proposition~\ref {pro:1.3} the exact
homotopy sequence of the fibration splits. Thus $f$ can be deformed
into the union of the fibre and the image of the section.
 \qed

\begin{Ex} \label {ex:2.2}
For $r=2r'>2$ and $m\ge1$ every pair of maps $f_1,f_2\:S^m\to
G_{r,2}$ (into the Grassmann manifold of $2$-planes in $\RR^r$) is
loose.

Indeed, due to the complex structure on $\RR^r=\CC^{r'}$, the
fibration $S^{r-2}\hookrightarrow V_{r,2}\to S^{r-1}$ (of the Stiefel
manifold of $2$-frames in $\RR^r$) has a canonical section; hence
 $$
 \pi_m(V_{r,2})\cong\pi_m(S^{r-2})\oplus\pi_m(S^{r-1}).
 $$
For $m\ge3$ this group is also isomorphic to $\pi_m(G_{r,2})$, and
the summands correspond to the subspaces
 $$
 A:=\{\mbox{$2$-planes containing the base vector $e_r$}\}
 $$
and
 $$
 B:=\{\mbox{complex lines in $\CC^{r'}$}\}.
 $$
Since $A\cup B\subsetneqq G_{r,2}$, our claim follows from lemma~\ref
{lem:2.1}.
 \qed
\end{Ex}

\proof{of proposition~\ref {pro:1.4}}
 A homotopy lifting argument shows that $(f,f)$ is loose by small
deformation if and only if the pulled back bundle $f^*(TN)$ has a
nowhere vanishing section over $S^m$ (compare~\cite{DG}). If $N$ is
noncompact or $\chi(N)=0$, then the tangent bundle $TN$ itself has a
nonzero section over $f(S^m)$ which we can pull back. In any case
every $n$-plane bundle over $S^m$ allows a section with a single
zero; it can be removed if its ``index map'' $q\:S^{m-1}\to S^{n-1}$
(compare~\cite{Wy} or also~\cite{K6}, (28)) is nullhomotopic.
 \qed

\medskip
The claim in example~\ref {ex:1.5} follows since $\pi_{m-1}(S^1)=0$
except when $m=2$. Furthermore $\pi_2(N)=0$ for every surface other
than $N=S^2$ or $\RP(2)$.

\section{Looseness obstructions} \label{sec:3}

In this section we recall the definitions of the various versions of
the $\omega$-invariants. (For more details we refer to~\cite{K6}; see
also~\cite{K3} and~\cite{K2}.) Moreover we compare the Nielsen
components of $\omega^\#$ in the root case (in~\ref {pro:3.12}). This
will be needed to establish our main results.

Unless specified otherwise we will assume that $m,n\ge2$. Fix
basepoints $x_0\in S^m$ and $y_0\in N$, and a local orientation of
$N$ at $y_0$. Throughout the remainder of this paper we consider
(continuous) maps $f_1,f_2,f,\ldots\:(S^m,x_0)\to(N,y_0)$. Our
notation will not distinguish between a constant map and its value.

After a suitable approximation we may assume that $(f_1,f_2)\:S^m\to
N\times N$ is smooth and transverse to the diagonal
 $$
 \Delta = \{(y,y)\in N\times N\mid y\in N\}.
 $$
Then the coincidence set
$$ \refstepcounter {Thm} \label {eq:3.1}
C=\C(f_1,f_2)=(f_1,f_2)^{-1}(\Delta)=\{x\in S^m\mid f_1(x)=f_2(x)\}
 \leqno (\theThm)
$$
is a smooth closed $n$-codimensional submanifold of $S^m$. The map
$(f_1,f_2)$ induces an isomorphism of (normal) vector bundles
$$ \refstepcounter {Thm} \label {eq:3.2}
 \nu(C,S^m)\cong(f_1,f_2)^*(\nu(\Delta,N\times N))\cong f_1^*(TN)|_C.
 \leqno (\theThm)
$$
Now pick a homotopy $G\:C\times I\to S^m$ which deforms the
inclusion map $g\:C\hookrightarrow S^m$ to the constant map with
value $x_0$. (Such a contraction exists and is unique up to homotopy
rel $(0,1)$ since $n\ge2$ and the complement of any point in $S^m$
allows ``linear'' homotopies.) Then $G$ induces a vector bundle
isomorphism from $f_1^*(TN)|_C=g^*(f_1^*(TN))$ to \ $C\times
T_{y_0}(N)$. Composing this with~(\ref {eq:3.2}) and using our
choice of a local orientation of $N$ at $y_0$, we get a framing
$$ \refstepcounter {Thm} \label {eq:3.3}
 \bar g^\#\:\nu(C,S^m)\stackrel{\cong}{\longrightarrow} C\times\RR^n.
 \leqno (\theThm)
$$
Furthermore we obtain the map
$$ \refstepcounter {Thm} \label {eq:3.4}
 \tilde g\:C \to \Omega(N,y_0)=:\Omega N
 \leqno (\theThm)
$$
which assigns the (concatenated) loop
$$ \refstepcounter {Thm} \label {eq:3.5}
 y_0=f_1(x_0) \stackrel{(f_1\circ G(x,-))^{-1}}{\longrightsquigarrow}
     f_1(x)=f_2(x) \stackrel{f_2\circ G(x,-)}{\longrightsquigarrow}
            f_2(x_0) = y_0
 \leqno (\theThm)
$$
to $x\in C$.

The resulting bordism class
$$ \refstepcounter {Thm} \label {eq:3.6}
 \omega^\#(f_1,f_2)=[C,\bar g^\#,\tilde g]
 \leqno (\theThm)
$$
of the framed compact submanifold $C\subset S^m$ (cf.~(\ref
{eq:3.1});~(\ref {eq:3.3})) together with the map $\tilde g$
(cf.~(\ref {eq:3.4})) depends only on the homotopy classes
$[f_i]\in\pi_m(N,y_0)$, $i=1,2$. Via the Pontryagin-Thom procedure,
$\omega^\#(f_1,f_2)$ can also be interpreted as an element in the
$m$-th homotopy group of the Thom space $S^n\wedge(\Omega N)^+$ of
the trivial $n$-plane bundle over the loop space $\Omega N$. (Here
``$+$'' stands for a disjointly added point. Note that the bordism
theories of submanifolds in $S^n$ and $\RR^n$ are equivalent in
codimensions $n\ge2$; thus it is not necessary that $f_1(x_0)\ne
f_2(x_0)$, as was assumed in~\cite{K6}.)

If we ignore the map $\tilde g$ we obtain the simpler invariant
$$ \refstepcounter {Thm} \label {eq:3.7}
 \underline{\omega}^\#(f_1,f_2)=[C,\bar g^\#] \in \pi_m(S^n).
 \leqno (\theThm)
$$
However, often this means a considerable loss of information. Indeed,
in general the loop space $\Omega N$ has a very rich topology and, in
particular, can be highly disconnected. Already its decomposition
into pathcomponents leads to the important ``Nielsen decomposition''
of coincidence sets, as follows.

Given a pathcomponent $A\in\pi_0(\Omega N)=\pi_1(N,y_0)$, restrict
your attention to the corresponding partial coincidence manifold
 $$
 C_A:=\C_A(f_1,f_2):=\tilde g^{-1}(A) \ \  \subset \ \ \C(f_1,f_2)
 $$
which is again a closed $n$-codimensional submanifold of $S^n$ and
endowed with the restricted framing $\bar g_A^\#=\bar g^\#|$ and the
map $\tilde g_A=\tilde g|\:C_A\to A\subset\Omega N$. This leads to
the invariants
$$ \refstepcounter {Thm} \label {eq:3.8}
 \omega_A^\#(f_1,f_2)=[C_A,\bar g_A^\#,\tilde g_A] \in \pi_m(S^n\wedge A^+)
 \leqno (\theThm)
$$
and
$$ \refstepcounter {Thm} \label {eq:3.9}
 \underline{\omega}_A^\#(f_1,f_2)=[C_A,\bar g_A^\#] \in \pi_m(S^n).
 \leqno (\theThm)
$$

\begin{Ex} \label {ex:3.10}
For all $[f]\in\pi_m(N,y_0)$ we have
$\underline{\omega}^\#(y_0,f)=\coll_*([f])$. Here
 $$
 \coll\: N\to N/(N\setminus\Int{B^n}) = B^n/\d B^n = S^n
 $$
is the map which collapses the complement of a small ball $B^n$ near
$y_0\in\d B^n$ to a point and preserves the local orientation of $N$.
 \qed
\end{Ex}

In the root case we will now compare the various Nielsen components
of $\omega^\#$.

Note that there are two canonical isomorphisms
$$ \refstepcounter {Thm} \label {eq:3.11}
 \rho_{A*}, \lambda_{A*}\: \pi_m(S^n\wedge A_0^+) \stackrel{\cong}{\longrightarrow}
 \pi_m(S^n\wedge A^+)
 \leqno (\theThm)
$$
(compare~(\ref {eq:3.8})) where $A_0$ denotes the pathcomponent of
$\Omega N$ containing the constant loop. They are induced by the
homotopy equivalences $\rho_A,\lambda_A\:A_0\to A$ which compose each
loop $\ell\in A_0$ to the right (or left, resp.) with a fixed loop
$\ell_A\in A$; e.g.\ $\rho_A(\ell)$ travels first along $\ell$ and
then along $\ell_A$.

To simplify notation we will write $\omega_0^\#$ for
$\omega_{A_0}^\#$.

\begin{Pro} \label {pro:3.12}
Given a map $f\:(S^m,x_0)\to(N,y_0)$ and $A\in\pi_1(N)$, we have
 $$
 \omega_A^\#(f,y_0)=\rho_{A*}\left(\omega_0^\#(f,y_0)\right)
 \quad \mbox{and} \quad
 \underline{\omega}_A^\#(f,y_0)=\underline{\omega}_0^\#(f,y_0)
 $$
as well as
 $$
 \omega_A^\#(y_0,f)=\iota_{A*}\circ\lambda_{A*}\left(\omega_0^\#(f,y_0)\right)
 \quad \mbox{and} \quad
 \underline{\omega}_A^\#(y_0,f)=\iota_{A*}\left(\underline{\omega}_0^\#(y_0,f)\right)
 $$
where $\iota_{A*}$ composes the framing with an orientation
preserving (or reversing) automorphism of $\RR^n$ according as $A$
is (or is not) orientation preserving, i.e.\ the tangent bundle of
$N$, when pulled by $\ell_A\:S^1\to N$, $[\ell_A]=A$, becomes
trivial or, equivalently, $A$ lies in the kernel of the composed
homomorphism
 $$
 w_1(N)\:\pi_1(N,y_0) \to H_1(N) \to \ZZ_2
 $$
which evaluates the first Stiefel-Whitney class of $N$.
\end{Pro}

\pproof{for some of the following arguments compare also the proof of
theorem~4.3 in~\cite{K6}}

In view of proposition~\ref {pro:1.3}~(iii) we need to consider only
the case when $\pi_1(N)$ is finite. Thus the fibre $p^{-1}(\{y_0\})$
in the universal covering space $p\:\tilde N\to N$ consists of points
$\{\tilde y_1,\dots,\tilde y_k\}$ which all lie in a suitable
embedded $n$-ball $V$ in $\tilde N$. Now lift $f$ to a map $\tilde
f\:S^m\to\tilde N$. After a suitable homotopy we may assume that
$\tilde f$ is smooth, with regular value $\tilde y_1$, and agrees on
a tubular neighbourhood $U\cong\tilde C\times V\subset S^m$ of
$\tilde C:=\tilde f^{-1}(\{\tilde y_1\})$ with the projection to
$V\subset\tilde N$. Then
 $$
 \C(f,y_0)=\C(y_0,f) = f^{-1}(\{y_0\}) = \tilde f^{-1}(\{\tilde y_1,\dots,\tilde y_k\})
 $$
consists of the ``parallel'' copies
 $$
 \tilde C_i:=\tilde C\times \{\tilde y_i\} \subset U, \quad i=1,\dots,k.
 $$
Note that the straight path $\gamma_{ij}$ joining $\tilde y_i$ to
$\tilde y_j$ in the ball $V$ yields an isotopy which moves $\tilde
C_i$ to $\tilde C_j$ in $U\subset S^m$, $1\le i,j\le k$. We can
compose it with a contraction $G|_{\tilde C_j\times I}$ of $\tilde
C_j$ in order to get the required contraction of $\tilde C_i$ in
$S^m$. This extra part of the contraction induces a concatenation of
the path in~(\ref {eq:3.5}) with the loop $p\circ\gamma_{ij}^{\mp1}$.
On the other hand the isotopy is compatible with the framings of
$\tilde C_i$ and $\tilde C_j$ if $p\circ\gamma_{ij}$ preserves the
local orientation of $N$ or in case we are dealing with
$\omega^\#(f,y_0)$ since then
 $$
 \nu(\tilde C_t,S^n)\cong\tilde C_t\times TV\cong f^*(TN)|_{\tilde
 C_t},
 $$
$0\le t\le 1$, throughout the isotopy. Thus let us consider
$\omega^\#(y_0,f)$. Here $\tilde C_i$ and $\tilde C_j$ are framed via
the orientations of $T_{\tilde y_i}\tilde N$ and $T_{\tilde
y_j}\tilde N$, resp., induced by $p$ from the given orientation of
$N$ at $y_0$. If $p\circ\gamma_{ij}$ reverses it, then the two
framings correspond to opposite orientations of $V$, and the isotopy
does not change this.
 \qed

\medskip
Next let us recall how the involution \ $\inv$ \ in~(\ref {eq:1.9})
is defined (cf.~\cite{K6}, p.~632). Given $[C,\bar g^\#,\tilde
g]\in\pi_m(S^n\wedge(\Omega N)^+)$, \ $\inv$ \ retains the
submanifold $C\subset S^m$ but changes its framing by
$(-1)^n\cdot\alpha$ where the vector bundle automorphism $\alpha$ of
the trivial bundle $C\times\RR^n$ is determined by $TN$ and the
homotopy $C\times I\to N$ which evaluates $\tilde g$
(cf.~\cite{K3},~3.1). In addition, we have to compose $\tilde g$ with
the selfmap of $\Omega N$ which inverts the direction of loops.

\begin {Ex} [$\mathbf{m=n\equiv0(2)}$, $N$ orientable]  \label {ex:3.13}
 Here the framings
(= coorientations) of the (isolated) points $x\in C$ remain
unchanged, but we have to travel backwards along each loop $\tilde
g(x)$.
\end{Ex}

\begin{War} In general the symmetry relation~(\ref {eq:1.9})
does not necessarily extend to $\underline{\omega}^\#$: the weaker
invariant $\underline{\omega}^\#(f_2,f_1)$ may depend on more than
just its weak counterpart $\underline{\omega}^\#(f_1,f_2)$.
\end{War}

A weaker version of our looseness obstruction $\omega^\#(f_1,f_2)$
is often much easier to handle and to compute. If we forget about
{\it embeddings} of coincidence manifolds into $S^m$ and if we keep
track only of stabilized normal bundles we obtain the invariants
$\tilde\omega_A(f_1,f_2)$, $A\in\pi_1(N)$, and
$$ \refstepcounter {Thm} \label {eq:3.14}
 \tilde\omega(f_1,f_2)\ \ \ \in\ \ \ \Omega_{m-n}^\fr(\Omega N)\ \cong\ \pi_{m+q}(S^{n+q}\wedge(\Omega N)^+),
 \quad q\gg0,
 \leqno (\theThm)
$$
which lie in the $(m-n)$-dimensional framed bordism group of the loop
space $\Omega N$ of $N$ (cf.~\cite{K3}). Similarly
$$ \refstepcounter {Thm} \label {eq:3.15}
 \omega(f_1,f_2):=\tilde{\underline\omega}(f_1,f_2), \ \
 \tilde{\underline\omega}_A(f_1,f_2)\ \ \ \in\ \ \ \Omega_{m-n}^\fr = \pi_{m-n}^S
 \leqno (\theThm)
$$
(compare~\cite{K2}, 1.4) are the stabilized versions of
$\underline{\omega}^\#(f_1,f_2)$ and
$\underline{\omega}_A^\#(f_1,f_2)$.

In general this stabilization procedure leads to a loss of
information; not so, however, in the dimension range $m<2n-2$ \ where
$\tilde\omega(f_1,f_2)$ is just as strong as $\omega^\#(f_1,f_2)$
and, in fact, is the {\it complete} looseness obstruction for the
pair $(f_1,f_2)$ (cf.~\cite{K6}, (16) and~\cite{K3}, theorem~1.10).

\section{Selfcoincidences} \label{sec:4}

In this section we prove theorem~\ref {thm:1.21} and further
vanishing results and discuss their consequences.

Given any map $f\:(S^m,x_0)\to(N,y_0)$ let us look at the coincidence
data of the pair $(f,f)$. Clearly the map $\tilde g$ (cf.~(\ref
{eq:3.4})) is canonically nullhomotopic since $f_1=f_2=f$ in~(\ref
{eq:3.5}). For small generic approximations of $(f,f)$ the partial
coincidence manifolds $C_A$ are empty (and hence nullbordant)
whenever $A\in\pi_1(N)$ is nontrivial.

We conclude that
$$ \refstepcounter {Thm} \label {eq:4.1}
 \omega^\#(f,f)=s_*(\underline\omega_0^\#(f,f))
 \leqno (\theThm)
$$
where
 $$
 s_*\:\pi_m(S^n) \to \pi_m(S^n\wedge(\Omega N)^+)
 $$
is induced by the inclusion of the constant loop into $\Omega N$
(recall also the notation of~\ref {pro:3.12}).

\medskip
\proof{of theorem~\ref {thm:1.21}}
 We base our argument on the identity
$$ \refstepcounter {Thm} \label {eq:4.2}
 \underline\omega^\#_A(f,f)=\underline\omega_A^\#(f,y_0) + \underline\omega_A^\#(y_0,f)
 \leqno (\theThm)
$$
(valid for all $A\in\pi_1(N)$; cf.~(\ref {eq:1.13})) and on
proposition~\ref {pro:3.12}. The assumptions in~\ref {thm:1.21} just
mean that the homomorphism $w_1\:\pi_1(N)\to\ZZ_2$ (cf.~\ref
{pro:3.12}) is not injective. In other words, there exists an
element $A$ of $\pi_1(N)$ which is both nontrivial and orientation
preserving. Then on the one hand $\underline\omega^\#_A(f,f)=0$. On
the other hand~\ref {pro:3.12} and~(\ref {eq:4.2}) combine to show
that $\underline\omega^\#_A(f,f)=\underline\omega^\#_0(f,f)$. Thus
in view of~(\ref {eq:4.1}) the full $\omega^\#$-invariant of the
pair $(f,f)$ -- and of every pair freely homotopic to it (cf.\
remark~\ref {rem:1.32}) -- vanishes.

To complete the proof of theorem~\ref {thm:1.21} note that
 $$
 \omega^\#(f,f) = \omega^\#(f,y_0) + \omega^\#(y_0,f) = (\id +
 \inv)(\deg^\#(f))
 $$
(cf.~(\ref {eq:1.8}),~(\ref {eq:1.9}) and~(\ref {eq:1.10})). Recall
also corollary~\ref {cor:1.16}.
 \qed

\medskip
Corollary~\ref {cor:1.25} follows now from~\ref {pro:1.3}, \ref
{pro:1.4}, (\ref {eq:1.11}), (\ref {eq:1.13}), \ref {thm:1.15} (to be
discussed below), \ref {thm:1.21}, \ref{pro:1.24} and the
compatibility of $\underline\omega(f,f)$ with covering maps
(compare~(\ref {eq:4.1})).

As for the conclusion of example~\ref {ex:1.22} see also~\ref
{ex:3.13}. Similarly, the statement in example~\ref {ex:1.23} is a
consequence of theorem~\ref {thm:1.21} and the following wellknown

\begin{Fac} [real Grassmann manifolds] \label {fac:4.3}
Given integers $0<k<r$, the manifold $G_{r,k}$ of all $k$-dimensional
linear subspaces of $\RR^r$ enjoys the following properties:
\begin{enumerate}
\item $\pi_1(G_{r,k})\cong\ZZ_2$ if $r>2$;
\item $G_{r,k}$ is orientable if and only if $r$ is even;
\item $\dim(G_{r,k})=k(r-k)$;
\item the Euler characteristic vanishes if $k\not\equiv r\equiv0(2)$
and equals the binomial coefficient
 $$
 \chi(G_{r,k})={[r/2]\choose[k/2]}>0
 $$
in all other cases.
\end{enumerate}
\end{Fac}

\begin{Que} \label {que:4.4}
What about $\omega^\#(f,f)$ for maps into $G_{r,k}$ or $\tilde
G_{r,k}$ when $r$ is odd and $k>2$?
\end{Que}

\proof{of theorem~\ref {thm:1.15}}
 A homotopy lifting argument shows that $(f,f)$ is loose by small deformation
precisely if the pulled back vector bundle $f^*(TN)$ has a nowhere
vanishing section over $S^m$ (cf.~\cite{DG} or~\cite{K7},~5.3) or,
equivalently, if $f$ lifts to $ST(N)$ (cf.~(\ref {eq:1.14})). Thus
(i)$\iff$(ii).

Recall also from~\cite{K7},~5.7 that $\pm E(\d[f])$ equals the
invariant $\underline\omega^\#(f,f)$ which is just as strong as
$\omega^\#(f,f)$ (cf.~(\ref {eq:4.1})).

Next compare the fibre homotopy sequence of $ST(N)$ and of the
configuration space $\tilde C_2(N)=N\times N\setminus\Delta$
(cf.~\cite{K7},~5.4). We see that (iii)$\implies$(i) provided the
induced homomorphism $j_*$ (cf.~\ref {pro:1.24}) is injective. This
is the case e.g.\ when $N=\RP(n)$.

Finally recall that (iv) implies (iii) when $N=S^n$; this is a
special case of~\cite{K6},~1.12.
 \qed

\medskip
Next we prove corollary~\ref {cor:1.16}. Clearly $E$ (cf.~(\ref
{eq:1.14})) is injective if $m=n$ or $m\le2n-3$ or $n=2$ (since
$\pi_{m-1}(S^1)=0$ for $m>2$). Moreover $\d(\pi_m(N))=0$ if $n$ is
odd since then the fibration $ST(N)\to N$ allows a section over each
compactum (compare~(\ref {eq:1.14})). Thus our corollary follows from

\begin {Pro} \label {pro:4.5}
Consider the suspension homomorphism
 $$
 E\:\pi_{m-1}(S^{n-1}) \to \pi_m(S^n)
 $$
and
 $$
 E^\infty\:\pi_{m-1}(S^{n-1}) \to \pi_{m-n}^S
 $$
and assume that $n$ is even. If $m\le n+3$ or if $m=n+4\ne8,10$ then
$E$ and $E^\infty$ are injective. If $m=8$ and $n=4$, then $E$ is
injective, but $E^\infty$ is not.
\end{Pro}

\Proof
 We will use Toda's tables~\cite{T}. Since
$E^\infty\:\pi_5(S^3)\cong\ZZ_2\to\pi_2^S\cong\ZZ_2$ is onto and
hence bijective, it remains only to study the (nonstable) cases where
$n=4$ and $m=7$ or $8$. The groups in the exact EHP-sequence
 $$
 \pi_7(S^3) \stackrel{E}{\longrightarrow}
 \pi_8(S^4) \stackrel{H}{\longrightarrow}
 \pi_8(S^7) \stackrel{P}{\longrightarrow}
 \pi_6(S^3) \stackrel{E}{\longrightarrow}
 \pi_7(S^4) \longrightarrow \ldots
 $$
(cf.~\cite{W}, XII, 2.3) have order $2,4,2$ and $12$, resp.; hence
both homomorphisms $E$ are injective here. Moreover the kernel of
 $$
 E^\infty\:\pi_7(S^4)\cong\ZZ\oplus\ZZ_{12} \to \pi_3^S\cong\ZZ_{24}
 $$
is generated by $[\iota_4,\iota_4]$ and has a trivial intersection
with $E(\pi_6(S^3))=\{0\}\oplus\ZZ_{12}$. Our claim follows for
$(m,n)=(7,4)$ and similarly for $(m,n)=(8,4)$ (since
$\pi_7(S^3)\ne0=\pi_4^S$).
 \qed

\begin{Cor} \label {cor:4.6}
There exists a map $f\:S^8\to\RP(4)$ such that $\omega^\#(f,f)\ne0$
but $\tilde\omega(f,f)=0$. Of course the corresponding liftings
$\tilde f\:S^8\to S^4$ have the same property.
\end{Cor}

\Proof The groups in the exact sequence
 $$
 \pi_8(ST(S^4)) \to
 \pi_8(S^4) \stackrel{\d}{\longrightarrow}
 \pi_7(S^3)
 $$
(cf.~(\ref {eq:1.14})) have order 2 (cf.~\cite{P}), 4 and 2, resp.;
thus there is an element $[\tilde f]\in\pi_8(S^4)$ such that
$\d([\tilde f])$ and hence $\underline\omega^\#(\tilde f,\tilde
f)=E(\d[\tilde f])$ is nontrivial, but $\omega(\tilde f,\tilde
f)=\tilde{\underline\omega}(\tilde f,\tilde f)\in\pi_4^S=0$ (compare
also~(\ref {eq:4.1})).
 \qed

\medskip
Next we want to explore the fact that it is often easier to compute
the stabilized versions of our $\omega$-invariants. Indeed they just
sum up the contributions of the partial coincidence manifolds $C_A$
(not registering their linkings); also we just multiply a given
bordism class by $-1$ if we compose the framing with a reflection of
$\RR^n$.

\begin {Pro} \label {pro:4.8}
Given $[f]\in\pi_m(N)$, the stabilized invariant $\tilde\omega(f,f)$
is determined by $\tilde{\underline\omega}(f,f)$, i.e.\ by
 $$
 \omega(f,f) = \chi(N) \, \omega (f,y_0)
 = \omega (f,y_0) + E^\infty (\coll_*([f]))
 = \pm E^\infty \circ \d ([f]).
 $$
(Here $\chi(N)$ denotes the Euler characteristic of $N$; for \
$\coll_*$ see~$\ref {ex:3.10}$.) If the stable suspension
homomorphism $E^\infty\:\pi_{m-1}(S^{n-1})\to\pi_{m-n}^S$ is
injective (or if $m\le n+3$) then the conditions {\rm (i)--(v)} in
theorem~$\ref {thm:1.15}$ are all equivalent to

\smallskip
{\rm (vi)} $\omega(f,f)=0$.
\end{Pro}

\Proof
 The first identity was already established in~\cite{K2},~2.2; the
corresponding claim for $\omega^\#(f,f)$ (see~\cite{K6},~5.1) is more
complicated and not so easy to use in calculations.

The second identity follows from~(\ref {eq:1.13}) and~\ref {ex:3.10}.
Furthermore note that $\omega(f,f)=\underline{\tilde\omega}(f,f)$ is
the stable suspension of $\underline\omega^\#(f,f)=\pm E(\d([f]))$
(cf.~\cite{K7},~5.7).

If $m\le n+3$ and $n$ is even, then $E^\infty$ is injective (cf.~\ref
{pro:4.5}). If $n$ is odd, all conditions (i)--(vi) hold
automatically.
 \qed

\begin{Pro} \label {pro:4.9}
Assume that $k:=\#\pi_1(N)\ge2$. Consider arbitrary
maps\\
$f,f_1,f_2\:(S^m,x_0)\to(N,y_0)$.

If $N$ is orientable, then
 $$
 \chi(N) \cdot \omega(f_1,f_2) \ \ \ \ \ \ \ \ =\ 0
 \eqno\mbox{\rm (cf.~(\ref {eq:3.14}))}
 $$
and in particular
 $$
 \chi(N) \cdot E^\infty(\coll_*([f]))\ =\ 0
 \eqno\mbox{\rm (cf.~\ref {ex:3.10}).}
 $$

\medskip
If $N$ is not orientable then \; $E^\infty(\coll_*([f]))\ =\ 0$ \;
and
 $$
\ \  (\chi(N)-1) \ \omega(f_1,f_2)\ \ =\ 0\ ;
 $$
if in addition $k>2$, then \; $\omega(f_1,f_2)=0$.
\end{Pro}

\Proof Note that $\omega(y_0,f)=E^\infty(\coll_*([f]))$ (cf.~\ref
{ex:3.10}). If $N$ is orientable, then
$$ \refstepcounter {Thm} \label {eq:4.10}
 0=\omega(f,f)=\omega(f,y_0) + \omega(y_0,f) = \chi(N) \cdot \omega(f,y_0)
 \leqno (\theThm)
$$
(cf.~\ref {thm:1.21}, (\ref {eq:1.13}) and~\ref {pro:4.8}). Thus
multiplication by the Euler characteristic annihilates also
$\omega(y_0,f)=-\omega(f,y_0)$ and hence $\omega(f_1,f_2)$ (by~(\ref
{eq:1.13})).

If $N$ is not orientable, then it follows from proposition~\ref
{pro:3.12} that
$$ \refstepcounter {Thm} \label {eq:4.11}
 \tilde{\underline\omega}_A(f,y_0)=\tilde{\underline\omega}_0(f,y_0)
 \quad \mbox {and} \quad
 \tilde{\underline\omega}_A(y_0,f)=\eps_A\cdot\tilde{\underline\omega}_0(y_0,f)
 \leqno (\theThm)
$$
for every $A\in\pi_1(N)$ where $\eps_A=+1$ or $-1$ according as $A$
is orientation preserving or reversing. Thus
$$ 
 \omega(y_0,f)=\sum_A \tilde{\underline\omega}_A(y_0,f)=0.
$$
Moreover
 $$
\tilde{\underline\omega}_A(f,f)
 = \tilde{\underline\omega}_A(f,y_0) + \tilde{\underline\omega}_A(y_0,f)
 =0
 $$
whenever $A\ne0$ (cf.~(\ref {eq:1.13}) and~(\ref {eq:4.1})). If in
addition $k>2$ then there exist both orientation preserving and
reversing $A\ne0$ so that
$\tilde{\underline\omega}_0(f,y_0)=-\eps_A\cdot\tilde{\underline\omega}_0(y_0,f)$
both for $\eps_A=+1$ and $\eps_A=-1$; hence $\omega(f,y_0)$, being an
even multiple of an element of order 2, vanishes, and so does
$\omega(f_1,f_2)$, again by~(\ref {eq:1.13}). If $\pi_1(N)$ consists
only of $0$ and $A\ne0$, then
$$ \refstepcounter {Thm} \label {eq:4.12}
 \tilde{\underline\omega}_0(f,y_0)=\tilde{\underline\omega}_A(f,y_0)
 =-\tilde{\underline\omega}_A(y_0,f)=\tilde{\underline\omega}_0(y_0,f)
 \leqno (\theThm)
$$
and
$$ \refstepcounter {Thm} \label {eq:4.13}
 \omega(f,y_0)=2\tilde{\underline\omega}_0(f,y_0)
 =\omega(f,f)=\chi(N)\,\omega(f,y_0)
 \leqno (\theThm)
$$
(cf.~\ref {pro:4.8}), so that multiplication by $\chi(N)-1$
annihilates $\omega(f,y_0)$ and hence also
$\omega(f_1,f_2)=\omega(f_1,y_0)$ (cf.~(\ref {eq:4.11})).
 \qed

\medskip
Given a map $f\:(S^m,x_0)\to(N,y_0)$, consider a lifting $\tilde
f\:(S^m,x_0)\to(\tilde N,\tilde y_0)$ to the universal covering space
$p\:\tilde N\to N$.

\begin{Lem} \label {lem:4.14}
We have
 $$
 \tilde{\underline\omega}_0(f,y_0) = \omega(\tilde f,\tilde y_0)
 \quad \mbox {and} \quad
 \tilde{\underline\omega}_0(y_0,f) = \omega(\tilde y_0,\tilde f).
 $$
\end{Lem}

\Proof In the Nielsen decomposition
 $$
 f^{-1}(\{y_0\}) = \tilde f^{-1}(p^{-1}(\{y_0\}))
 = \tilde f^{-1}(\{A\tilde y_0\mid A\in\pi_1(N)\})
 $$
of the relevant coincidence manifold, $\tilde f^{-1}(\{\tilde y_0\})$
is the component indexed by $A=0$.
 \qed

\medskip
Let us apply this to the special case $\tilde N=S^n$ and $N=\RP(n)$,
$n$ even. Then
 $$
 \omega(f,f) = 2 \tilde{\underline\omega}_0(y_0,f)
 = 2 \omega(\tilde y_0,\tilde f) = 2E^\infty([\tilde f])
 $$
(cf.~(\ref {eq:4.13}), \ref {lem:4.14}, \ref {ex:3.10} and~(\ref
{eq:3.15})). This establishes the claim which follows theorem~\ref
{thm:1.31}. Also if $m=2n-1$ and e.g.\ $n=4,8,12,14,16$ or $20$, then
$\ker E^\infty\cong\ZZ$ and $E^\infty$ maps $\pi_m(S^n)$ onto the
stable stem $\pi_{n-1}^S$ which contains elements of order greater
than $2$ (cf.~\cite{T}). This proves the statement in example~\ref
{ex:1.26}.

Similarly, (\ref {eq:4.13}) shows also (together with theorem~\ref
{thm:1.21}) that $\omega(f,f)\in2\cdot\pi_{m-n}^S$ whenever $N$ is
not simply connected. This (together with~\ref {pro:4.9})
establishes~(\ref {eq:1.27}).

Finally let us discuss example~\ref {ex:1.28}. According to
fact~\ref {fac:4.3}, the Grassmann manifold $G_{5,2}$ is not
orientable, $6$-dimensional and its Euler characteristic equals~$2$.
Our claim follows from proposition~\ref {pro:4.9} and the
Freudenthal suspension theorem.

\section{Nielsen numbers} \label{sec:5}

\begin{Def} \label {def:5.1}
Given $f_1,f_2\:S^m\to N$, the {\it ``strong'' Nielsen number}
$\N^\#(f_1,f_2)$ (and its stabilized analogue \ $\N(f_1,f_2)$,
resp.) is the number of elements $A\in\pi_1(N)$ such that
$\omega^\#_A(f_1,f_2)$ (and $\tilde\omega_A(f_1,f_2)$, resp.) does
not vanish (cf.~(\ref {eq:3.8}) and~(\ref {eq:3.14})).
\end{Def}

Since generic coincidence manifolds are compact these Nielsen numbers
are always finite.

\medskip
\proof{of theorem~\ref {thm:1.30}} If $\omega^\#(f_1,f_2)$ (or
$\tilde\omega(f_1,f_2)$, resp.) vanishes, then so do all the partial
invariants $\omega^\#_A(f_1,f_2)$ (or $\tilde\omega_A(f_1,f_2)$,
resp.), $A\in\pi_1(N)$, as well as the corresponding Nielsen number.
The converse is also obvious in the stabilized setting. However,
$\omega^\#(f_1,f_2)$ keeps track also of embeddings and of possible
linking phenomena among the partial coincidence submanifolds
$\C_A(f_1,f_2)$ in $S^m$.

Assume that $\N^\#(f_1,f_2)=0$. Then all triples
$(\C_A(f_1,f_2),\bar g_A^\#,\tilde g_A)$ (cf.~(\ref {eq:3.8})) admit
individual nullbordisms in $S^m\times I$. Conceivably these can not
be fitted together disjointly to yield the full {\it embedded}
nullbordism required to show that $\omega^\#(f_1,f_2)$ vanishes. In
fact, it is an open question whether the first claim in~\ref
{thm:1.30} still holds when the common domain of $f_1$ and $f_2$ is
not a sphere.

But here we can use the additive structure of homotopy groups. As
in~(\ref {eq:1.8}) we have for all $A\in\pi_1(N)$
$$ \refstepcounter {Thm} \label {eq:5.2}
 0=\omega^\#_A(f_1,f_2)=\omega^\#_A(f_1,f_1) + \omega^\#_A(y_0,f_2-f_1).
 \leqno (\theThm)
$$
The claim of our theorem is obvious when $\pi_1(N)=0$ or in the
selfcoincidence case since only $A=0$ plays a r\^ole there (cf.~(\ref
{eq:4.1})). If $A\ne0$ then $\omega^\#_A(f_1,f_1)=0$ and hence
$\omega^\#_A(y_0,f_2-f_1)=0$ (cf.~(\ref {eq:5.2})). This implies --
in this special root case -- that the full invariant
$\omega^\#(y_0,f_2-f_1)$ vanishes; indeed, a nullbordism of
$\C_A(y_0,f_2-f_1)$ gives rise to {\it disjoint} ``parallel''
nullbordisms of all the other partial coincidence manifolds
$\C_{A'}(y_0,f_2-f_1)$, $A'\in\pi_1(N)$ (compare the proof of~\ref
{pro:3.12} or also~\cite{K6},~4.3). Thus in turn
$\omega^\#(f_1,f_2)=\omega^\#(f_1,f_1)$ (cf.~(\ref {eq:1.8})) and we
are back in the selfcoincidence case.
 \qed

\medskip
\proof{of theorem~\ref {thm:1.31}} If $k=\#\pi_1(N)$ is infinite,
then all pairs $(f_1,f_2)$ are loose and hence all Nielsen numbers
vanish (cf.~(\ref {eq:3.1})). If $k>2$ or if $k=2$ and $N$ is
orientable, then always $\omega^\#(f_1,f_1)=0$ (cf.~\ref {thm:1.21});
therefore $\omega^\#(f_1,f_2)=\omega^\#(y_0,f_2-f_1)$ (cf.~(\ref
{eq:5.2})) and $\tilde\omega(f_1,f_2)=\tilde\omega(y_0,f_2-f_1)$.
Again~\ref {pro:3.12} implies that the corresponding Nielsen numbers
can assume only the values $0$ and $k$ in this root case. For the
remaining cases compare corollary~\ref {cor:1.25}.
 \qed

\begin{Rem} \label {rem:5.3}
Consider the split exact sequence
 $$
 0\to\ker(\forg)\hookrightarrow\pi_m(S^n\wedge(\Omega N)^+)
 \stackrel{\forg}{\longrightarrow} \bigoplus_{A\in\pi_1(N)}\pi_m(S^n\wedge A^+)\to0\ \ ;
 $$
here, given any element $\omega^\#=[C,\bar g^\#,\tilde
g]\in\pi_m(S^n\wedge(\Omega N)^+)$, $\forg(\omega^\#)$ keeps track
only of the individual $A$-components $[C_A=\tilde g^{-1}(A),\bar
g^\#|,\tilde g|]$ and forgets about possible linkings. As we will
see below, the kernel of \ $\forg$ \ can be highly nontrivial if
$\#\pi_1(N)>1$.
\end{Rem}

Consider also the function $\N^\#$ which counts the nontrivial
$A$-components. Clearly it vanishes on $\ker(\forg)$ and, in
addition, can often assume all integer values between $0$ and
$\#\pi_1(N)$. This shows that theorems~\ref {thm:1.30} and~\ref
{thm:1.31} impose strong restrictions on those $\omega^\#$-values
which can actually be realized by pairs $(f_1,f_2)$ of maps.

\begin{Ex} [real projective spaces] \label {ex:5.4}
Consider $N=\RP(n)$ and its double cover $\tilde N=S^n$. There is a
wellknown isomorphism
 $$
 \pi_m(S^n\wedge(\Omega N)^+)\cong\pi_m(S_1^n\vee S_2^n\vee\tilde N,\tilde N)
 $$
(cf.\ e.g.~\cite{K6},~(61)), and the forgetful map \ $\forg$ \
(cf.~\ref {rem:5.3}) corresponds to the homomorphism
 $$
 \pi_m(\bigvee^2S^n\vee\tilde N) \to \bigoplus_{i=1}^2\pi_m(S_i^n\vee\tilde N)
 $$
which is induced by the obvious projections. Let $\iota_1,\iota_2$
and $\iota_0$ be represented by the inclusions of the three spheres
into their wedge. Then the summands of
 $$
 \ker(\forg)
 \cong \pi_m(S^{2n-1})\oplus(\pi_m(S^{3n-2}))^4\oplus\ldots
 $$
in the Hilton decomposition of $\pi_m(S_1^n\vee S_2^n\vee\tilde N)$
(cf.~\cite{W},~XI,~6) corresponds precisely to those basic products
which involve both $\iota_1$ and $\iota_2$.
\end{Ex}

\end{document}